\documentclass[12pt]{article}
\usepackage{amsmath,amssymb,amsthm}
\usepackage[margin=1.25in]{geometry}
\usepackage[bookmarksnumbered]{hyperref}
\usepackage{enumitem,theoremref,verbatim,cancel, breqn, color, authblk}

\title{Projections of Gibbs States for H\"older Potentials}
\author{Mark Piraino}
\affil{University of Victoria}

\theoremstyle{definition}
\newtheorem{theorem}{Theorem}

\newtheorem{lemma}[theorem]{Lemma}
\newtheorem{example}[theorem]{Example}

\newtheorem{proposition}[theorem]{Proposition}
\newtheorem{definition}[theorem]{Definition}
\newtheorem*{notation}{Notation}
\newtheorem*{remark}{Remark}
\newtheorem*{acknowledgments}{Acknowledgments}

\newcommand{\bmat}[1]{\begin{bmatrix} #1 \end{bmatrix}}
\newcommand{\inn}[1]{\left\langle #1 \right\rangle}
\newcommand{\set}[1]{\left\{ #1 \right\}}
\newcommand{\abs}[1]{\left| #1 \right|}
\newcommand{\norm}[1]{\left \| #1 \right \|}

\def\[#1\]{\begin{align*}#1\end{align*}}
\DeclareMathOperator{\var}{var}
\DeclareMathOperator{\diam}{diam}

\newcommand{\N}{\mathbb{N}}

\newcommand{\R}{\mathbb{R}}

\def\B{\mathcal{B}}     
\def\A{\mathcal{A}}      
\def\X{\mathcal{X}}  
\def\L{\mathcal{L}}

\begin{document}
	\maketitle
	
\begin{abstract}
	In this paper we give a short proof that the projection of a Gibbs state for a H\"older continuous potential on a mixing shift of finite type under a 1-block fiber-wise mixing factor map has a H\"older continuous g function. This improves a number of previous results. The key insight in the proof is to realize the measure of a cylinder set in terms of positive operators and use cone techniques. 
\end{abstract}

\section{Introduction}

On shifts of finite type one distinguished class of measures are Gibbs states for H\"older continuous potentials (H\"older Gibbs states). These measures have their origins in statistical physics and play an important role in the ergodic theory of axiom A diffeomorphisms. For an introduction to these measures we refer the reader to Bowen's famous monograph \cite{bowen1975equilibrium}. Somewhat surprisingly the projection of a H\"older Gibbs state under a continuous factor map need not be a H\"older Gibbs state. We will give an example of this later. A natural question is then to ask what conditions on the factor map ensure that H\"older Gibbs states project to H\"older Gibbs states?

Incremental progress has been made in many directions, for a $1$-block code $\pi$ between full shifts and $\varphi$ H\"older it was shown in \cite{chazottes2009preservation} and \cite{pollicott2011factors} that $\pi_{\ast}\mu_\varphi$ is the Gibbs state for a function $\psi$ with $\var_{n}(\psi)=O(\eta^{\sqrt{n}})$ where $0<\eta <1$. In fact in \cite{pollicott2011factors} the more general case of potential of summable variation was considered. This was improved in \cite{zbMATH05938049} to give that $\var_{n}(\psi)=O(\eta^{n})$ answering our question for factors of the full shift. 

On shifts of finite type this question was first posed in \cite{chazottes2003projection}, where the case of Markov measures was considered. It has been shown that under some conditions on $\pi$, the projection of Markov measures are H\"older Gibbs states (\cite{chazottes2003projection},\cite{Yoo2010}). The work in \cite{pollicott2011factors} was generalized to shifts of finite type in \cite{kempton2011factors} again under some assumptions on the factor map giving that $\pi_{\ast}\mu_{\varphi}$ is a Gibbs state for a function $\psi$ with $\var_{n}(\psi)=O(\eta^{\sqrt{n}})$. We significantly improve this result and settle the problem by showing that the projection of any H\"older Gibbs state under a fiber-wise mixing factor map (definition \ref{fiberwisemixingdef}) is a H\"older Gibbs state. The insight that allows our proof to work is that the measure of a cylinder set for a H\"older Gibbs state can be computed in terms of positive operators acting on $C(\Sigma_{A} \to \R)$. The theory of hidden Markov (sofic) measures can then be adapted to our setting, in particular our approach could be viewed as the analog of \cite{chazottes2003projection} with infinite dimensional spaces.

Throughout the article let $\A,\B$ be finite alphabets, $\Sigma_{A}\subseteq \A^{\N}$ a topologically mixing shift of finite type, $\varphi:\Sigma_{A}\to \R$ such that $\var_{n}\varphi\leq \abs{\varphi}_{\theta}\theta^{n}$ and $\pi:\A \to \B$ a map inducing a 1-block factor $\pi:\Sigma_{A} \to Y \subseteq \B^{\N}$.

\begin{definition}\label{fiberwisemixingdef}
	Call a $1$-block factor map, $\pi$, fiber-wise mixing if there exists an $N$ such that for any admissible word $b_{0}\cdots b_{N}$ in $Y$ and $a_{0},a_{N}\in \A$ such that $\pi(a_{0})=b_{0}$ and $\pi(a_{N})=b_{N}$, there exists a word $a_{0}a_{1}\cdots a_{n-1}a_{N}$ admissible in $\Sigma_{A}$ with $\pi(a_{0}\cdots a_{N})=b_{0}\cdots b_{N}$.
\end{definition}

Without the condition of being fiber-wise mixing there are examples of Markov Measures which don't project to H\"older Gibbs states \cite{chazottes2003projection}. We will give our own example \ref{exampleofMarkovtoNonHolder} below. Projections under fiber-wise mixing factors are however well behaved. In particular we have the following theorem.

\begin{theorem}\label{mainthm}
	Suppose that $\varphi:\Sigma_{A} \to \R$ is H\"older continuous, $\mu_{\varphi}$ the Gibbs state for $\varphi$, and $\pi:\Sigma_{A} \to Y$ a fiber-wise mixing $1$-block factor map. Then the projected measure $\pi_{\ast}\mu_{\varphi}$ is the Gibbs state for a H\"older continuous potential.
\end{theorem}

One advantage of our approach is that for the full shift we have a formula for the regularity of a potential associated to the projected measures in terms of the regularity of the original potential. That is if $\var_{n}\varphi \leq \abs{\varphi}_{\theta}\theta^{n}$ then $\pi_{\ast}\mu_{\varphi}$ is a Gibbs state for a function $\psi$ with $\var_{n}(\psi)=O(\eta^{n})$ where
\[ \eta = \tanh\left(\frac{1}{2}\left( \log \frac{1+\sigma}{1-\sigma}+\sigma \frac{\abs{\varphi}_{\theta}\theta}{\sigma - \theta}\right)\right)  \]
for any $\theta < \sigma <1$. We begin by fixing some notation.

\begin{notation}
	Write
	\[ X_{i}=\set{f \in C(\Sigma_{A} \to \R):f(x)=0 \text{ for all }x \in \Sigma_{A}\setminus [i]}. \]
	This is a subspace of $C(\Sigma_{A} \to \R)$, moreover $C(\Sigma_{A} \to \R)=\bigoplus_{i \in \A}X_{i}$. For $b \in \B$ define
	\[ \X_{b}=\bigoplus_{\pi(i)=b}X_{i} =\set{f \in C(\Sigma_{A} \to \R):f(x)=0 \text{ for all }x \in \Sigma_{A}\setminus \bigcup_{\pi(i)=b} [i]}. \]
	For $f \in \X_{b}$ we write $f_{i}$ for the component in $X_{i}$.
\end{notation}

For $i,j \in \A$ with $A_{ij}=1$ define $L_{ij}:X_{i} \to X_{j}$
\[ L_{ij}f(x)=e^{\varphi(ix)}f(ix)\chi_{[j]}(x), \]
where $ix$ is the point defined by $(ix)_{k}=i$ if $k=0$ and $(ix)_{k}=x_{k-1}$ otherwise. Setting $L_{ij}=0$ otherwise and identifying $C(\Sigma_{A} \to \R) = \bigoplus_{i\in \A} X_{i}$ we can see that the transfer operator can be written as
\[ L_{\varphi}=\bmat{L_{11}&L_{21}&\cdots &L_{n1}\\ L_{12}&L_{22}&\cdots &L_{n2}\\ \vdots &\vdots & &  \\ L_{1n}& L_{2n}& &L_{nn} }. \]
Take $h$ and $\nu$ such that $L_{\varphi}h=\lambda h$ and $L^{\ast}_{\varphi}\nu = \lambda \nu$ and $\inn{h,\nu}=1$. 

\begin{lemma}\label{GibbsStateOpRep}
	Suppose that $\varphi:\Sigma_{A}\to \R$ is H\"older and $\mu_{\varphi}$ is the Gibbs state. Then
	\[ \mu_{\varphi}([x_{0}\cdots x_{n}])= \lambda^{-n}\inn{L_{x_{n-1}x_{n}}\cdots L_{x_{0}x_{1}}h_{x_{0}},\nu_{x_{n}}}.  \]
	\begin{proof}
		Notice
		\begin{align*}
		\mu_{\varphi}[x_{0}\cdots x_{n}]=\lambda^{-n}\int L^{n}_{\varphi}(h\chi_{[x_{0}\cdots x_{n}]})d\nu&=\lambda^{-n}\int \sum_{\abs{I}=n,Ix\in \Sigma_{A}}e^{S_{n}\varphi(Iz)}h(Iz)\chi_{[x_{0}\cdots x_{n-1}x_{n}]}(Iz)d\nu(z)\\
		&=\lambda^{-n}\int_{[x_{n}]} e^{S_{n}\varphi(x_{0}\cdots x_{n-1}z)}h(x_{0}\cdots x_{n-1}z)d\nu(z).
		\end{align*}
		On the other hand
		\[L_{x_{1}x_{2}}L_{x_{0}x_{1}}h(z)&= e^{\varphi(x_{1}z)}L_{x_{0}x_{1}}h(x_{1}z)\chi_{[x_{2}]}(z)\\
		&=e^{\varphi(x_{1}z)}e^{\varphi(x_{0}x_{1}z)}h(x_{0}x_{1}z)\chi_{[x_{1}]}(x_{1}z)\chi_{[x_{2}]}(z)\\
		&=e^{\varphi(x_{1}z)}e^{\varphi(x_{0}x_{1}z)}h(x_{0}x_{1}z)\chi_{[x_{2}]}(z). \]
		By iteration we can see
		\[ L_{x_{n-1}x_{n}}\cdots L_{x_{1}x_{2}}L_{x_{0}x_{1}}h(z)= e^{S_{n}\varphi(x_{0}\cdots x_{n-1}z)}h(x_{0}\cdots x_{n-1}z)\chi_{[x_{n}]} \]
		and $\nu_{x_{n}}$ is the functional $\inn{f,\nu_{x_{n}}}= \int_{[x_{n}]}f d \nu$.
	\end{proof}
\end{lemma}

This gives us a way of describing the projected measure. For $b,b' \in \B$ with $bb'$ admissible in $Y$ we define $\L_{bb'}:\X_{b}\to \X_{b'}$ by
\[ \L_{bb'}f=\sum_{\pi(i)=b,\pi(j)=b'}L_{ij}f_{i}  \]
and $h_{b}=\sum_{\pi(i)=b}h_{i}$, $\nu_{b} = \sum_{\pi(i)=b}\nu_{i}$. The key observation is the following. For any word $b_{0}\cdots b_{n}\in  \L(Y)$
\[ \pi_{\ast}\mu_{\varphi}[b_{0}\cdots b_{n}]=\lambda^{-n}\inn{\L_{b_{n-1}b_{n}}\L_{b_{n-2}b_{n-1}}\cdots \L_{b_{0}b_{1}}h_{b_{0}},\nu_{b_{n}}}. \]
To see this simply notice that 
\begin{align*}
\L_{b_{n-1}b_{n}}\L_{b_{n-2}b_{n-1}}\cdots \L_{b_{0}b_{1}}h_{b_{0}}= \sum_{\pi(a_{0}\cdots a_{n})=b_{0} \cdots b_{n}}L_{a_{n-1}a_{n}}L_{a_{n-2}a_{n-1}}\cdots L_{a_{0}a_{1}}h_{a_{0}}.
\end{align*}
Thus
\begin{align*}
\inn{\L_{b_{n-1}b_{n}}\L_{b_{n-2}b_{n-1}}\cdots \L_{b_{0}b_{1}}h_{b_{0}},\nu_{b_{n}}}= \sum_{\pi(i)=b_{n},\pi(j)=b_{n}}\sum_{\pi(a_{0}\cdots a_{n-1})=b_{0}\cdots b_{n-1}}\inn{L_{a_{n-1}j}L_{a_{n-2}a_{n-1}}\cdots L_{a_{0}a_{1}}h_{a_{0}},\nu_{i}},
\end{align*}
where of course the duality pairing is $0$ unless $i=j$. Finally by lemma \ref{GibbsStateOpRep} we have
\[ \lambda^{-n}\sum_{\pi(a_{0}\cdots a_{n})=b_{0}\cdots b_{n}}\inn{L_{a_{n-1}a_{n}}L_{a_{n-2}a_{n-1}}\cdots L_{a_{0}a_{1}}h_{a_{0}},\nu_{a_{n}}}=\sum_{\pi(a_{0}\cdots a_{n})=b_{0}\cdots b_{n}}\mu_{\varphi}[a_{0}\cdots a_{n}]=\pi_{\ast}\mu_{\varphi}[b_{0}\cdots b_{n}]. \]
To prove that the projection of a H\"older Gibbs state is a H\"older Gibbs state we need a candidate for the potential on $Y$. The obvious choice is the $g$ function for $\pi_{\ast}\mu_{\varphi}$ which is given by the formula
\[ g(x)=\lim_{n \to \infty}\frac{\pi_{\ast}\mu_{\varphi}[x_{0}x_{1}\cdots x_{n}]}{\pi_{\ast}\mu_{\varphi}[x_{1}\cdots x_{n}]}=\lambda^{-1}\lim_{n \to \infty}\frac{\inn{\L_{x_{n-1}x_{n}}\cdots \L_{x_{1}x_{2}}\L_{x_{0}x_{1}}h_{x_{0}}, \nu_{x_{n}}}}{\inn{\L_{x_{n-1}x_{n}}\cdots \L_{x_{1}x_{2}}h_{x_{1}}, \nu_{x_{n}}}} \]
where it is understood that the function may only be defined $\pi_{\ast}\mu_{\varphi}$ almost everywhere. If the $g$ function is regular enough then $\pi_{\ast}\mu_{\varphi}$ is the Gibbs state associated to $\log g$. For more information on the connection between $g$ functions and equilibrium states see \cite{MR1085356}. Therefore it is enough for us to show that $\log g$ is H\"older. One can see that the regularity of $\log g$ is intimately connected to the convergence of 
\[ \frac{\L_{b_{1}b_{2}}^{\ast}\cdots \L_{b_{n-1}b_{n}}^{\ast}\nu_{b_{n}}}{\inn{h_{b_{1}},\L_{b_{1}b_{2}}^{\ast}\cdots \L_{b_{n-1}b_{n}}^{\ast}\nu_{b_{n}}}}. \]
This being a projective limit leads naturally to the use of cone techniques to tackle the problem. We illustrate the idea with an example of a factor map which is not fiber-wise mixing for which the projection of Parry measure does not have a H\"older $g$ function.

\begin{example}\label{exampleofMarkovtoNonHolder}
Let $\Sigma_{A}$ be the shift determined by the matrix
\[ A=\bmat{1&1&1&0\\ 0&1&1&1\\ 1&1&1&0\\ 0&1&1&1} \]
define $\pi:\Sigma_{A}\to \Sigma_{2}$ by $\pi(0)=\pi(1)=0$ and $\pi(2)=\pi(3)=1$, and take $\varphi = 0$. Consider the word $0^{k+1}\in \L(\Sigma_{2})$ and notice that if $\pi(a_{0}\cdots a_{k})=0^{k+1}$ and $a_{0}=1$ then $a_{i}=1$ for all $0 \leq i \leq k$ because $10\notin \L(\Sigma_{A})$. Therefore $\pi$ is not fiber-wise mixing. 

Consider the projection of the Parry measure on $\Sigma_{A}$ under $\pi$. In the case of Markov measures the transfer operator preserves a finite dimensional subspace of $C(\Sigma_{A}\to \R)$ containing the leading eigenvector, the functions locally constant on cylinder sets of length $1$. Therefore in this case we can work with matrices. It is easy to see that $\lambda =3$, $h=(1,1,1,1)$, $\nu=6^{-1}(1,2,2,1)$,
\[ \L_{00}=\bmat{1&1\\0&1}, \L_{01}=\bmat{1&1\\0&1}, \L_{10}=\bmat{1&0\\1&1}, \L_{11}=\bmat{1&0\\1&1} \]
and
\[ h_{0}=\bmat{1\\1} \text{ , }h_{1}=\bmat{1\\1} \text{ and }\nu_{0}=6^{-1}\bmat{1\\2} \text{ , }\nu_{1}=6^{-1}\bmat{2\\1}. \]
We can compute
\[ g(0^{\infty})=\frac{1}{3}\lim_{n\to \infty}\frac{\inn{\bmat{1&1\\0&1}^{n}\bmat{1\\1}, \bmat{1\\2}}}{\inn{\bmat{1&1\\0&1}^{n-1}\bmat{1\\1}, \bmat{1\\2}}}=\frac{1}{3}\lim_{n\to \infty}\frac{n+3}{n+2}=\frac{1}{3} \]
and taking $x_{k}=0\cdots 01111\cdots$ the point of $k-1$ zeros and then $1$'s again we compute $g(x_{k})$
\[ g(x_{k})&=\frac{1}{3}\lim_{n \to \infty}\frac{\inn{\bmat{1&0\\1&1}^{n-k}\bmat{1&1\\0&1}^{k}\bmat{1\\1},\bmat{2\\1}}}{\inn{\bmat{1&0\\1&1}^{n-k-1}\bmat{1&1\\0&1}^{k-1}\bmat{1\\1},\bmat{2\\1}}}\\
&=\frac{1}{3}\lim_{n \to\infty}\frac{n(k+1)+k-k^{2}+3}{nk+k-k^{2}+1}=\frac{1}{3}\left(\frac{k+1}{k}\right). \]
Therefore
\[ \abs{g(0^{\infty})-g(x_{k})}=\frac{1}{3k} \]
and $g$ cannot be H\"older.
\end{example}

This shows us that even in the case when $\varphi =0$ the projection of $\mu_{\varphi}$ under a continuous factor map can fail to be a H\"older Gibbs state. It also begs an interesting question, which is, can the assumption that $\pi$ is fiber-wise mixing be replaced with the assumption that the projection of the Parry measure is a H\"older Gibbs state in theorem \ref{mainthm}? Furthermore is it the case that $\pi$ is fiber-wise mixing if and only if the projection of the Parry measure is a H\"older Gibbs state? It seems very likely that the answer to the first question is yes. To answer the second, one needs to understand under what conditions a measure defined by products of non-negative matrices has a H\"older $g$ function. Next we review some properties of Hilbert's projective metric which we will need to prove theorem \ref{mainthm}.

\begin{definition}
	Let $V$ be a real Banach space. A subset $\Lambda \subseteq V$ is called a cone if 
	\begin{enumerate}
		\item 
		$\Lambda \cap (-\Lambda) = \set{0}$
		
		\item
		$c \Lambda \subseteq \Lambda$ for all $c \geq 0$
		
		\item
		$\Lambda$ is convex
	\end{enumerate}
	In addition we define the dual cone
	\[ \Lambda^{\ast}=\set{\phi \in V^{\ast}:\inn{x,\varphi}\geq 0 \text{ for all }x \in \Lambda} \]
\end{definition}

On any cone there is a notion of a projective distance called the Hilbert metric which we now define.

\begin{definition}
	Let $\Lambda$ be a closed cone. For $x,y \in \Lambda$ define
	\[ \alpha(x,y)=\sup\set{\lambda>0: y-\lambda x \in \Lambda} \text{ and }\beta(x,y)=\inf\set{\lambda>0:\lambda x - y \in \Lambda} \]
	where $\alpha(x,y)=0$ and $\beta(x,y)=\infty$ if the sets are empty. The Hilbert metric on $\Lambda$ is defined by
	\[ \Theta_{\Lambda}(x,y)=\log \left(\frac{\beta(x,y)}{\alpha(x,y)}\right). \]
\end{definition}
This is a projective pseudo-metric in the sense that it has the properties of a metric when restricted to the unit sphere, although it can take the value $\infty$, and for any $x,y\in \Lambda$ and $a,b>0$ we have $\Theta_{\Lambda}(ax,by)=\Theta_{\Lambda}(x,y)$. The true utility of this metric is the following famous theorem.

\begin{theorem}\label{BirkhoffContraction}
	(Birkhoff \cite{BirkhoffHilbertMetric}) Let $\Lambda_{1},\Lambda_{2}$ be closed cones and $L:V_{1} \to V_{2}$ a linear map such that $L\Lambda_{1} \subseteq \Lambda_{2}$. Then for all $\phi , \psi \in \Lambda_{1}$
	\[ \Theta_{\Lambda_{2}}(L\phi , L\psi) \leq \tanh\left(\frac{\diam_{\Lambda_{2}}(L\Lambda_{1})}{4}\right)\Theta_{\Lambda_{1}}(\phi , \psi) \]
	where 
	\[ \diam_{\Lambda_{2}}(L\Lambda_{1})=\sup\set{\Theta_{\Lambda_{2}}(f,g ):f,g\in L\Lambda_{1}} \]
	and $\tanh \infty =1$.
\end{theorem}

Birkhoff's contraction theorem has a long history of use in dynamics, mostly for proving upper bounds on the rate of convergence of the transfer operator and decay of correlations (\cite{MR2083432}, \cite{LiveraniDOC}, \cite{kondah1997vitesse}). The underlying principle in these papers is simple: find a cone whose image under the transfer operator has finite diameter. We encounter an additional subtlety which is that our operators will map between different cones but this is not a serious technical issue and the same ideas apply. We will need the following proposition to relate the projective metric to the variations of $\log g$.

\begin{proposition}\label{DualconeMetricComputation}
	Let $\Lambda$ be a closed cone and $x,y \in \Lambda$ such that $\Theta_{\Lambda}(x,y)<\infty$. Then for any $\phi \in \Lambda^{\ast}$, $\inn{x,\phi}=0$ if and only if $\inn{y,\phi}=0$ and
	\[ \Theta_\Lambda(x,y) = \log \left(\sup\set{\frac{\inn{x, \phi}\inn{y,\psi}}{\inn{y,\phi}\inn{x,\psi}}:\psi,\phi \in \Lambda^{\ast}\text{ and }\inn{y,\phi}\inn{x,\psi}\neq 0}\right)\]
	\begin{proof}
	The proof can be found in \cite{eveson_nussbaum_1995} lemma 1.4.
	\end{proof}
\end{proposition}

\section{Proof of Theorem \ref{mainthm}}
\begin{notation}
	The proof relies on the use of the following cones
	\[ \X_{b}^{+}=\set{f \in \X_{b}:f \geq 0} \]
	and
	\[ \Lambda_{b}^{K}=\set{f \in \X_{b}^{+}:f(x)\leq f(y)e^{K\theta^{n}} \text{ whenever }n\geq 1 \text{ and } x_{i}= y_{i}\text{ for all }0 \leq i \leq n-1} \]
\end{notation}

\begin{lemma}\label{boundonHilbertmetricSFTHolder}
	\begin{enumerate}
		\item
		Let $0<\sigma <1$ and suppose that $f,g \in \Lambda_{b}^{\sigma K}$ then 
		\[ \Theta_{\Lambda_{b}^{K}}(f,g)\leq 2 \log \left(\frac{1+\sigma}{1-\sigma}\right)+\Theta_{\X_{b}^{+}}(f,g) \]
		
		\item
		Let $N$ be as in definition \ref{fiberwisemixingdef} take $K =\frac{\abs{\varphi}_{\theta}}{\sigma - \theta^{N}}\sum_{i=1}^{N}\theta^{i} $ where $0<\theta^{N} <\sigma<1$, set $\Lambda_{b}:=\Lambda_{b}^{K}$. There is a constant $M$ such that for any admissible word $b_{0}\cdots b_{N}$ in $Y$
		\[\diam_{\Lambda_{b_{N}}}(\L_{b_{N-1}b_{N}}\cdots \L_{b_{0}b_{1}}\Lambda_{b_{0}})\leq M < \infty\]
	\end{enumerate}
	\begin{proof}
		\begin{enumerate}
			\item 
			This is contained in \cite{MR2083432} proposition 5.3.
			
			\item
			 Suppose that $f,g \in \Lambda_{b_{0}}$ and denote $\L:=\L_{b_{N-1}b_{N}}\cdots \L_{b_{0}b_{1}}$. If $x_{i}= y_{i}$ for $0\leq i \leq k-1$, $k\geq 1$ then
			\begin{align*}
			\L f(x)&=\sum_{\abs{I}=N, \pi(I)=b_{0}\cdots b_{N-1}}e^{S_{N}\varphi(Ix)}f(Ix)\\
			&=\sum_{\abs{I}=N, \pi(I)=b_{0}\cdots b_{N-1}}e^{S_{N}\varphi(Ix)-S_{N}\varphi(Iy)+S_{N}\varphi(Iy)}\frac{f(Ix)}{f(Iy)}f(Iy)\\
			&\leq \exp\left[\abs{\varphi}_{\theta}\sum_{i=1}^{N}\theta^{i+k}+K\theta^{N+k}\right]\L f(y)\\
			&= \exp\left[\sigma K \theta^{k}\right]\L f(y).
			\end{align*}
			Thus $\L f,\L g \in \Lambda_{b_{N}}^{\sigma K}$ and it remains to estimate $\Theta_{X_{b_{N}}^{+}}(\L f,\L g)$. This is where we use in a crucial way that the factor map is fiber-wise mixing. Choose a point $z \in \bigcup_{\pi(j)=b_{0}}[j]$ such that $f(z)= \norm{f}$. By the fiber-wise mixing condition, for any $i$ with $\pi(i)=b_{N}$ there is a word $w_{z}^{i}:=z_{0}a_{1}\cdots a_{N}$ with $a_{N}= i$ and $\pi(w_{z}^{i})=b_{0}\cdots b_{N}$. Thus for any $x \in \bigcup_{\pi(i)=b_{N}}[i]$
			\begin{align*}
			\L f(x)=\sum_{\abs{I}=N, \pi(I)=b_{0}\cdots b_{N-1}}e^{S_{N}\varphi(Ix)}f(Ix) \geq e^{-N\norm{\varphi}_{\infty}}f(w_{z}^{x_{0}}x)\geq e^{-N\norm{\varphi}_{\infty}-\theta K}\norm{f}.
			\end{align*}
			Therefore we have
			\begin{align*}
			\Theta_{X_{b_{N}}^{+}}(\L f,\L g)&=\log \left(\sup_{x,y \in \bigcup_{\pi(j)=b_{0}}[i]}\frac{\L f(x)\L g(y)}{\L g(x)\L f(y)}\right)\\
			&\leq \log \left(\frac{\norm{\L}_{op
			}^{2}\norm{f}\norm{g}}{e^{-2N\norm{\varphi}_{\infty}-2\theta K}\norm{f}\norm{g}}\right)\\
			&= \log \left(\frac{\norm{\L}_{op}^{2}}{e^{-2N\norm{\varphi}_{\infty}-2\theta K}}\right)\\
			&\leq \log \left(\frac{\norm{L_{\varphi}^{N}1}_{\infty}^{2}}{e^{-2N\norm{\varphi}_{\infty}-2\theta K}}\right)<\infty.
			\end{align*}
			The bounds are independent of the word $b_{0}\cdots b_{N}$, hence the result.
		\end{enumerate}
	\end{proof}
\end{lemma}

\begin{remark}
To get the estimate on $\eta$ for the full shift notice that any $1$-block factor map is fiber-wise mixing with $N=1$. The first estimate in lemma \ref{boundonHilbertmetricSFTHolder} holds for $k=0$ that is
\[\L f(x) \leq \exp\left[\sigma K \right]\L f(y)\]
for any $x,y \in \Sigma$. Hence
\[ \Theta_{\X_{b}^{+}}(\L f(x), \L f(y)) \leq 2 \sigma K  \]
\end{remark}

We are now in a position to prove the main theorem. Lemma \ref{boundonHilbertmetricSFTHolder} tells us that long enough products of our operators are contractions. What remains is to relate the variations of $\log g$ to the projective metric and apply Birkhoff's contraction theorem.

\begin{proof}
	(of Theorem \ref{mainthm}) Let $n \geq 0$ and take $x,y \in Y$ such that $x_{i}=y_{i}$ for $0 \leq i \leq n-1$. Consider for $m,k \geq n$
	\begin{align*}
	&\abs{\log \frac{ \inn{\L_{x_{m-1}x_{m}}\cdots \L_{x_{1}x_{2}}\L_{x_{0}x_{1}}h_{x_{0}}, \nu_{x_{m}}}}{\inn{\L_{x_{m-1}x_{m}}\cdots \L_{x_{1}x_{2}}h_{x_{1}}, \nu_{x_{m}}}}-\log \frac{\inn{\L_{y_{k-1}y_{k}}\cdots \L_{y_{1}y_{2}}\L_{y_{0}y_{1}}h_{y_{0}}, \nu_{y_{k}}}}{\inn{\L_{y_{k-1}y_{k}}\cdots \L_{y_{1}y_{2}}h_{y_{1}}, \nu_{y_{k}}}}}\\
	&=\abs{\log \frac{\inn{\L_{x_{n-1}x_{n}}\cdots \L_{x_{1}x_{2}}\L_{x_{0}x_{1}}h_{x_{0}}, \mu_{x,m}}\inn{\L_{x_{n-1}x_{n}}\cdots \L_{x_{1}x_{2}}h_{x_{1}},\mu_{y,k}}}{\inn{\L_{x_{n-1}x_{n}}\cdots \L_{x_{1}x_{2}}h_{x_{1}},\mu_{x,m}}\inn{\L_{x_{n-1}x_{n}}\cdots\L_{x_{1}x_{2}}\L_{x_{0}x_{1}}h_{x_{0}}, \mu_{y,k}}}}
	\end{align*}
	where
	\[ \mu_{x,m}=\L_{x_{n}x_{n+1}}^{\ast}\cdots \L_{x_{m-1}x_{m}}^{\ast}\nu_{x_{m}} \text{ and }\mu_{y,k}=\L_{y_{n}y_{n+1}}^{\ast}\cdots \L_{y_{k-1}y_{k}}^{\ast}\nu_{y_{k}} \]
	Notice that $\mu_{x,m},\mu_{y,k}\in \Lambda_{x_{n}}^{\ast}$. Assuming that $n\geq N$ as in definition \ref{fiberwisemixingdef} so that Hilbert's metric is finite we have by proposition \ref{DualconeMetricComputation} that
	\begin{align*}
	&\abs{\log \frac{\inn{\L_{x_{n-1}x_{n}}\cdots \L_{x_{1}x_{2}}\L_{x_{0}x_{1}}h_{x_{0}}, \mu_{x,m}}\inn{\L_{x_{n-1}x_{n}}\cdots \L_{x_{1}x_{2}}h_{x_{1}},\mu_{y,k}}}{\inn{\L_{x_{n-1}x_{n}}\cdots \L_{x_{1}x_{2}}h_{x_{1}},\mu_{x,m}}\inn{\L_{x_{n-1}x_{n}}\cdots\L_{x_{1}x_{2}}\L_{x_{0}x_{1}}h_{x_{0}}, \mu_{y,k}}}}\\
	&\leq \log \left(\sup\set{\frac{\inn{\L_{x_{n-1}x_{n}}\cdots \L_{x_{1}x_{2}}\L_{x_{0}x_{1}}h_{x_{0}}, \phi}\inn{\L_{x_{n-1}x_{n}}\cdots \L_{x_{1}x_{2}}h_{x_{1}},\psi}}{\inn{\L_{x_{n-1}x_{n}}\cdots \L_{x_{1}x_{2}}h_{x_{1}},\phi}\inn{\L_{x_{n-1}x_{n}}\cdots \L_{x_{1}x_{2}}\L_{x_{0}x_{1}}h_{x_{0}},\psi}}
	}\right)\\
	&=\Theta_{\Lambda_{x_{n}}}(\L_{x_{n-1}x_{n}}\cdots \L_{x_{1}x_{2}}\L_{x_{0}x_{1}}h_{x_{0}},\L_{x_{n-1}x_{n}}\cdots \L_{x_{1}x_{2}}h_{x_{1}})
	\end{align*}
	Where the supremum in the second line is over $\phi,\psi\in  \Lambda_{x_{n}}^{\ast}$ such that 
	\[\inn{\L_{x_{n-1}x_{n}}\cdots \L_{x_{1}x_{2}}h_{x_{1}},\phi}\inn{\L_{x_{n-1}x_{n}}\cdots \L_{x_{1}x_{2}}\L_{x_{0}x_{1}}h_{x_{0}},\psi} \neq 0. \]
	Write $n-1=qN+r$ with $0 \leq r \leq N-1$ and set $\eta=\tanh(M/4)$ where $M$ is as in lemma \ref{boundonHilbertmetricSFTHolder}. Birkhoff's contraction theorem gives us that
	\begin{align*}
	\Theta_{\Lambda_{x_{n}}}&(\L_{x_{n-1}x_{n}}\cdots \L_{x_{1}x_{2}}\L_{x_{0}x_{1}}h_{x_{0}},\L_{x_{n-1}x_{n}}\cdots \L_{x_{1}x_{2}}h_{x_{1}})\\
	&\leq \eta^{q-1}\Theta_{\Lambda_{x_{N+1}}}(\L_{x_{N}x_{N+1}}\cdots \L_{x_{1}x_{2}}\L_{x_{0}x_{1}}h_{x_{0}},\L_{x_{N}x_{N+1}}\cdots \L_{x_{1}x_{2}}h_{x_{1}}) \\
	&\leq \eta^{q-1}M= (\eta^{1/N})^{n} M\eta^{-1-(r+1)/N} \leq  (\eta^{1/N})^{n} M\eta^{-2}.
	\end{align*}
	Taking $x=y$ we see that the sequence 
	\[ \log \frac{\inn{\L_{x_{n-1}x_{n}}\cdots \L_{x_{1}x_{2}}\L_{x_{0}x_{1}}h_{x_{0}}, \nu_{x_{n}}}}{\inn{\L_{x_{n-1}x_{n}}\cdots \L_{x_{1}b_{2}}h_{x_{1}}, \nu_{x_{n}}}} \]
	is Cauchy and therefore $\log g(x)$ exists at every point. Moreover taking $m=k$ and letting $m \to \infty$ we have $\var_{n}(\log g)=O(\eta^{n/N})$ hence the result.
\end{proof}

Finally we briefly mention a strategy to extend theorem \ref{mainthm} beyond H\"older potentials. Cone techniques have been used to prove sub-exponential upper bounds on the rate of convergence for transfer operators associated to Walters functions \cite{kondah1997vitesse}. It seems reasonable to adapt our arguments and those of \cite{kondah1997vitesse} to obtain estimates on the regularity of $g$ functions associated to projections of Gibbs states for Walters functions. We leave this for future work.

\begin{acknowledgments}
I am grateful to Anthony Quas and Chris Bose for many useful discussions, as well as careful readings and helpful suggestions on several drafts of this manuscript.
\end{acknowledgments}

\bibliographystyle{abbrv}
\bibliography{FactorsHGS}{}

\end{document}